\documentclass[10pt]{article}
\usepackage{amsfonts}
\usepackage{mathrsfs}
\usepackage{graphicx, subfigure}
\usepackage{amsmath}
\usepackage{amssymb}
\usepackage[mathscr]{eucal}
\renewcommand{\paragraph}{\roman{paragraph}}
\newtheorem{theorem}{\scshape \bf Theorem}[section]
\newtheorem{lemma}[theorem]{\scshape \bf Lemma}

\begin{document}

\title{\LARGE \sf On the determinant of the distance matrix of a bicyclic graph\thanks
{Supported by National Natural Science Foundation of China(11071002,11171373),
  and Zhejiang Provincial Natural
Science Foundation of China(LY12A01016).}}
 \author{Shi-Cai Gong\thanks{Corresponding author. E-mail addresses:
scgong@zafu.edu.cn(S. Gong); zhangjl@zafu.edu.cn(L. Zhang); ghxu@
zafu.edu.cn(G. Xu).}, Ju-Li Zhang and Guang-Hui
Xu\vspace{3mm} \\
  School of Science, Zhejiang A \& F University,\\
 Lin'an, 311300, P. R. China
   }
\date{}
\maketitle

\noindent{\bf Abstract:} Two cycles   are referred as disjoint if
they have no common edges. In this paper, we will investigate the
determinant of the distance matrix of a  graph, giving a formula for
the determinant of the distance matrix of a bicyclic graph whose two
cycles are disjoint, which extends the formula for the determinant
of the distance matrix of   a tree, as well as that of  a unicyclic
graph. \vspace{2mm}

 \noindent
{\bf Keywords:} Distance matrix; bicyclic graph;  determinant

\noindent {\bf AMS Subject Classifications:} 05C50; 15A18

\section{Introduction} In the whole paper all graphs are simple and
undirected. Let $G$ be a graph with vertex set $V=\{1,2,\cdots,n\}$
and edge set $E$. The distance between the vertices $i$ and $j$,
denoted by $dis(i, j )$, is the length of a shortest path between
them. The $n$-by-$n$ matrix $D(G)=(d_{i,j})$ with $d_{i,j}=dis(i, j
)$ is referred as the  {\it distance matrix}  of $G$, or the metrics
matrix of $G$.

The determinant of the distance matrices of graphs have been
investigated in the literature. As early,  Graham and Pollack
\cite{gp} showed that if $T$ is a tree on $n$ vertices with distance
matrix $D$, then the determinant of $D$ is $(-1)^{n-1}(n-1)2^{n-2}$,
a formula depending only on $n$. Then Bapat, Kirkland and Neumann
\cite{bk} extend this formula to the weighted case and give a
formula for the determinant of the distance matrix of a unicyclic
graph, showing that the determinant of the distance matrix of a
unicyclic graph
 is related to the length of the cycle contained in it and its
 order. For
more spectral properties of the distance matrices of a graph, one
can see for example \cite{l,merris,y1,y2,zg}
 and the references therein.

For a given graph $G$, two cycles of $G$ are referred as disjoint if
they have no common edges. Let $C_p$ and $C_q$ be two disjoint
cycles. Suppose that $v_1\in C_p, v_k\in C_q $. Joining $v_1$ and
$v_k$ by a path $v_1v_2\cdots v_k$ of length $k-1$, where $k\ge 1$
and $k = 1$ means identifying $v_1$ with $v_k$, the resultant graph,
denoted by $\infty(p, k, q)$, is referred as an $\infty$-graph. A
bicyclic graph which contains an $\infty$-graph as an induced
subgraph can be considered a graph obtained from an $\infty$-graph
$\infty(p, k, q)$ by planting some trees to such an $\infty$-graph.

In this paper, we will investigate the determinant of the distance
matrix of a  graph, giving a formula for the determinant of the
distance matrix of a bicyclic graph whose two cycles are disjoint,
which extends the formula for the determinant of the distance matrix
of   a tree, as well as that of  a unicyclic graph. In addition, as
by-product we show  that if a graph is obtained from an induced
subgraph by planting some trees on it, then the determinant of the
distance matrix of such a graph is independent to the structure of
those trees.

\section{Preliminary results}
In this section, we will establish some preliminary results, which
will be useful in the following discussion.

 Henceforth we
use the following notation. For a real matrix $A$, denote by $A^T$
the transpose matrix to $A$.  The identity matrix is denoted by $I$
and the all ones row vector is denoted by $\textbf{1}$. The
determinant of the matrix $A$ is denoted by $det(A)$, or $|A|$ for
simplify.  We refer to D. Cvetkovi$\acute{c}$, M. Doob and H. Sachs
\cite{cds} for more terminology and notation not defined here.

\begin{lemma} \label{a0} Let $C_k=\frac{1}{2}B_kB_k^T-2I$, a $k\times k$ matrix,
 and $F_k=\frac{1}{2}\textbf{1}B^T_k+\textbf{1}$, a row vector with dimension
 $k$,
where {\small  $$B_k=\left(\begin{array}{cccccc}
-1&0&0&\cdots&0&0\\
-1&-1&0&\cdots&0&0\\
0&-1&-1&\cdots&0&0\\
\cdots& \cdots& \cdots& \cdots& \cdots& \cdots\\
0&0&0&\cdots&-1&0\\
0&0&0&\cdots&-1&-1
\end{array}\right)_{k\times k}.$$} Then $$det C_k=\frac{(-1)^k(2k+1)}{2^k},$$ and
$$F_kC_k^{-1}F_k^T=-\frac{k}{2(2k+1)}.$$
\end{lemma}
{\bf Proof.} By a directly calculation, we have {\small $$
C_k=\frac{1}{2}\left(\begin{array}{cccccc}
-3&1& 0&\cdots&0&0\\
1&-2&1&\cdots&0&0\\
0& 1& -2&\cdots&
0&0\\
\cdots&\cdots&\cdots&\cdots&\cdots&\cdots\\
0&0&0&\cdots&-2&1\\
0&0&0&\cdots&1&-2\\
\end{array}\right)_{k\times k}$$ and $$F_k=\frac{1}{2}(1~0~\cdots
~0~0),$$ a vector with exactly one nonzero entry. Now let
\begin{eqnarray*} H_k=\left(\begin{array}{cccccc}
-2&1& 0&\cdots&0&0\\
1&-2&1&\cdots&0&0\\
0& 1& -2&\cdots&
0&0\\
\cdots&\cdots&\cdots&\cdots&\cdots&\cdots\\
0&0&0&\cdots&-2&1\\
0&0&0&\cdots&1&-2\\
\end{array}\right)_{k\times k}.
\end{eqnarray*}}
As we know that $det H_k=(-1)^k(k+1),$ then {\small
\begin{eqnarray*} 2^k det C_k&=&\left|\begin{array}{cccccc}
-3&1& 0&\cdots&0&0\\
1&-2&1&\cdots&0&0\\
0& 1& -2&\cdots&
0&0\\
\cdots&\cdots&\cdots&\cdots&\cdots&\cdots\\
0&0&0&\cdots&-2&1\\
0&0&0&\cdots&1&-2\\
\end{array}\right|\\&=&\left|\begin{array}{cccccc}
-2&1& 0&\cdots&0&0\\
1&-2&1&\cdots&0&0\\
0& 1& -2&\cdots&
0&0\\
\cdots&\cdots&\cdots&\cdots&\cdots&\cdots\\
0&0&0&\cdots&-2&1\\
0&0&0&\cdots&1&-2\\
\end{array}\right|+\left|\begin{array}{cccccc}
-1&0& 0&\cdots&0&0\\
1&-2&1&\cdots&0&0\\
0& 1& -2&\cdots&
0&0\\
\cdots&\cdots&\cdots&\cdots&\cdots&\cdots\\
0&0&0&\cdots&-2&1\\
0&0&0&\cdots&1&-2\\
\end{array}\right|\\&=&det H_k-det H_{k-1}\\&=&(2k+1)(-1)^k.
\end{eqnarray*}}
Hence, $$det C_k=\frac{(2k+1)(-1)^k }{2^k},$$ and
 $$F_kC_k^{-1}F_k^T=\frac{det(C_{1,1}^*)}{4|C_k|}=\frac{det(H_{k-1})}{4|C_k|}=-\frac{k}{2(2k+1)},$$where
 $C_{1,1}^*$ denotes the $(1,1)-$th entry of the adjoint matrix of $C_k$. The result thus
 follows. \hfill
$\blacksquare$\vspace{3mm}

\begin{lemma} \label{a00} Suppose that the sequence $f(0),f(1),\cdots, f(n)$ satisfies the following linear recurrence
relation$$\left \{ \begin{array}{l} f(n)=-4f(n-1)-4f(n-2)\\
f(0)=f_0\\
f(1)=f_1.
\end{array}\right.$$Then $$f(n)=[f_0-\frac{n}{2}(f_1+2f_0)](-2)^n.$$
\end{lemma}
{\bf Proof.} Since the characteristic equation of this recurrence
relation is
 $$x^2+4x+4=0$$ and its two roots are $x_1=x_2=-2$,  by Theorem
 7 .4.1 in \cite{b} the general solution is
$$f(n)=(c_1+nc_2)(-2)^n.$$ Combining with the initial values
$ f(0)=f_0$ and $ f(1)=f_1$, we have $$\left \{ \begin{array}{l}
c_1=f_0;\\
c_2=-\frac{1}{2}f_1-f_0.
\end{array}\right.$$ The result thus follows.\hfill
$\blacksquare$\vspace{3mm}

\begin{lemma} \label{a1} Let $G$  be the  graph  obtained from
 a graph $G_1$  by identifying an arbitrary vertex of $G_1$
and one  pendent vertex of the path $P_2$.
 Then the determinant
of the distance matrix of the graph $G$ is fixed, regardless the
choice of the vertex of $G_1$.
\end{lemma}
{\bf Proof.} Let  the vertex set of $G$ be $\{1,2,\cdots,n\}$.
Without loss of generality, we can take vertex $1$ to be one pendent
vertex of $P_2$ and label another pendent vertex, a quasi-pendent
vertex of $G$, as $2$. Then $G_1$ can be considered the subgraph of
$G$ induced by vertices $\{2,\cdots,n\}$. Let $(0~~d_2)$ be the row
vector of the distance matrix   of $G_1$ corresponding to the vertex
$2$ and $D^*$ be the distance matrix of the subgraph of $G$ induced
by vertices $\{3,\cdots,n\}$. Then $D(G)$, the distance matrix   of
$G$, can be partitioned as
{\small  $$D(G)=\left(\begin{array}{ccc} 0& 1& d_2+\textbf{1}\\
1&0& d_2 \\
d_2^T+\textbf{1}^T& d_2^T & D^*
\end{array}\right).$$} Hence
{\small  $$det D(G)=\left|\begin{array}{ccc} 0& 1& d_2+\textbf{1}\\
1&0& d_2 \\
d_2^T+\textbf{1}^T& d_2^T & D^*
\end{array}\right|=\left|\begin{array}{ccc}  -1& 1& \textbf{1}\\
1&0& d_2 \\
d_2^T+\textbf{1}^T& d_2^T & D^*
\end{array}\right|=\left|\begin{array}{ccc}  -2& 1& \textbf{1}\\
1&0& d_2 \\
\textbf{1}^T& d_2^T & D^*
\end{array}\right|,$$}the last equalition implies that $det D(G)$  is independent to
the choice of the vertex $2$.  The result thus follows. \hfill
$\blacksquare$\vspace{3mm}

\begin{lemma} \label{a3} Let $G_1$ and $G_2$ be two graphs with vertex sets
 $\{1,2,\cdots,k\}$ and
$\{k+1,k+2\cdots,n\}$, respectively.
 Let $G$ be the graph
obtained from $G_1$ and $G_2$ by adding an edge between  vertices
 $1$ and $n$, and $\tilde{G}$  the graph
obtained from $G_1$ and $G_2$ by identifying vertices $1$ and $n$
and then adding a pendent vertex from  $1$ {\em (or $n$)}. Denote by
$D$ and $\tilde{D}$ respectively the distance matrices of $G$ and
$\tilde{G}$. Then $$det D=det \tilde{D}.$$
\end{lemma}
{\bf Proof.} Without loss of generality, we take the distance
matrices of $G_1$ and $G_2$
as  $D(G_1)=\left(\begin{array}{cc} 0&d_1\\
d_1^T&D^*
\end{array}\right)$ and
$D(G_2)=\left(\begin{array}{cc} D^{**}&d_n^T\\
d_n&0
\end{array}\right)$, where $D^*$ and $D^{**}$ denote respectively
the distance matrix of the subgraphs induced by $\{2,\cdots,k\}$ and
$\{k+1,k+2\cdots,n-1\}$, and $(0 ~~d_1)$ and $(d_n~~0)$ are
respectively the row vectors of $D(G_1)$ corresponding to the vertex
$1$ and the row vectors of $D(G_2)$ corresponding to the vertex $n$.
Again without loss of generality, suppose that, in $\tilde{G}$, the
vertex $n$ is the pendent vertex and the vertex $1$ is the
quasi-pendent vertex. For $D=(d_{i,j})$ and
$\tilde{D}=(\tilde{d}_{i,j})$, we set the rows and columns of them
correspond to $\{1,2,\cdots,n\}$, respectively. Then we have {\small
$$d_{i,j}=\left \{ \begin{array}{ll}
dis_{G_1}(i,j), & {\rm if \mbox{ } 1\le i,j\le k};\\
dis_{G_2}(i,j),& {\rm if \mbox{ } k+1\le i,j\le
n};\\
dis_{G_1}(1,i)+dis_{G_2}(j,n)+1,& {\rm if \mbox{ } 1\le i\le k
\mbox{ } k+1\le j\le n}, \end{array}\right.$$and
$$\tilde{d}_{i,j}=\left \{ \begin{array}{ll}
dis_{G_1}(i,j), & {\rm if \mbox{ } 1\le i,j\le k};\\
dis_{G_2}(i,j),& {\rm if \mbox{ } k+1\le i,j\le
n-1};\\
dis_{G_1}(1,j)+1,& {\rm if \mbox{ } i=n \mbox{ } and \mbox{
} 1\le j\le k};\\
dis_{G_2}(n,j)+1,& {\rm if \mbox{ } i=n \mbox{ } and \mbox{
} k+1\le j\le n-1};\\
dis_{G_1}(1,i)+dis_{G_2}(j,n),& {\rm if \mbox{ } 1\le i\le k \mbox{
} and \mbox{ } k+1\le j\le n-1}. \end{array}\right.$$}
  Hence,
{\small $$D=\left(\begin{array}{cccc}
0&d_1&d_n+\textbf{1}& 1\\
d_1^T& D^* &d_n\textbf{1}^T+\textbf{1}d_1^T+\textbf{1}\textbf{1}^T& d_1^T+\textbf{1}^T \\
d_n^T+\textbf{1}^T& \textbf{1}d_n^T+d_1\textbf{1}^T+\textbf{1}\textbf{1}^T &D^{**}& d_n^T\\
1& d_1+\textbf{1}& d_n&0
\end{array}\right),$$
 $$\tilde{D}=\left(\begin{array}{cccc}
0&d_1&d_n& 1\\
d_1^T& D^* &d_n\textbf{1}^T+\textbf{1}d_1^T& d_1^T+\textbf{1}^T\\
d_n^T& \textbf{1}d_n^T+d_1\textbf{1}^T &D^{**}& d_n^T+\textbf{1}^T\\
1& d_1+\textbf{1} & d_n+\textbf{1}&0
\end{array}\right).\eqno{(2.1)}$$} Consequently, we have
{\small \begin{eqnarray*} det D&=&\left|\begin{array}{cccc}
0&d_1&d_n+\textbf{1}& 1\\
d_1^T& D^* &d_n\textbf{1}^T+\textbf{1}d_1^T+\textbf{1}\textbf{1}^T& d_1^T+\textbf{1}^T \\
d_n^T+\textbf{1}^T& \textbf{1}d_n^T+d_1\textbf{1}^T+\textbf{1}\textbf{1}^T &D^{**}& d_n^T\\
1& d_1+\textbf{1}& d_n&0
\end{array}\right|\\&=&\left|\begin{array}{cccc}
0&d_1&d_n& 1\\
d_1^T& D^*-\textbf{1}d_1^T-d_1\textbf{1}^T &0& d_1^T\\
d_n^T& 0 &D^{**}-\textbf{1}d_n^T-d_n\textbf{1}^T& d_n^T\\
1& d_1 & d_n&0
\end{array}\right|\\&=&\left|\begin{array}{cccc}
0&d_1&d_n& 1\\
d_1^T& D^* &d_n\textbf{1}^T+\textbf{1}d_1^T& d_1^T+\textbf{1}^T\\
d_n^T& \textbf{1}d_n^T+d_1\textbf{1}^T &D^{**}& d_n^T+\textbf{1}^T\\
1& d_1+\textbf{1} & d_n+\textbf{1}&0
\end{array}\right|=det \tilde{D}.
\end{eqnarray*}} Then the result follows.\hfill
$\blacksquare$\vspace{3mm}

\section{  On the determinant of the distance matrix of a bicyclic graph whose two cycles are disjoint}

For a bicyclic graph $G$, if its two cycles are disjoint, then $G$
contains  $\infty(p,k,q)$ as an induced subgraph for some integers
$p$, $q$ and $k$.  This subgraph $\infty(p,k,q)$ is sometimes called
the center construct of $G$. In this way, the graph $G$ can be
viewed as the graph obtained from $\infty(p,k,q)$ by planting some
trees on it. In the following discussion, the  graph $\infty(p,1,q)$
will play an important role. For convenience,  the vertex, in
$\infty(p,1,q)$, with degree  $4$ is called the center vertex of
$\infty(p,1,q)$, and denote by $G(p,q;n)$ the graph obtained from
$\infty(p,1,q)$ by planting the path $P_{n}$ on its center vertex.
Then $G(p,q;0)$ denotes $\infty(p,1,q)$ itself and the graph
$G(p,q;n)$ has order $n+p+q-1$.\vspace{3mm}

First of all, combining with Lemmas \ref{a1} and \ref{a3}, we have
the following result, which tell us   that if a graph is obtained
from an induced subgraph by planting some trees on it, then the
determinant of the distance matrix of such a graph is independent to
the structure of those trees.
\begin{theorem} \label{a4} Let $G$  be a bicyclic graph of order $n+p+q-1$
which contains $\infty(p,k,q)$  as an induced subgraph for some
integers $p$, $q$ and $k$. Suppose that $D$ and $\tilde{D}$ be
respectively the distance matrices of $G$ and $G(p,q;n)$.
 Then
$$det D=det \tilde{D}.$$
\end{theorem}
{\bf Proof.} First, applying Lemma \ref{a3} repeatedly, the distance
matrices corresponding respectively to the graphs $\infty(p,k,q)$
and $G(p,q;k-1)$ have the same determinant.

Then it remain to show that the bicyclic graph, denoted still by
$G$, of order $n+p+q-1$ which contains $\infty(p,1,q)$  as an
induced subgraph has the same determinant as that of the graph
$G(p,q;n)$. We label the vertices of $G$ as $\{1,2,\cdots,p+q+n-1\}$
such that the resultant graph obtained from $G$ by deleting the
vertices $\{n,n-1,\cdots,n-i\}$ with $i(0\le i\le n)$  is connected.
We first consider the vertex $n$, if $n$ is not adjacent to the
center vertex of $G$, then applying Lemma \ref{a1} to $G$ such that
the vertex $n$ adjacent to the center vertex of $G$, the resultant
graph is still denoted by $G$; then applying Lemma \ref{a1} to $G$
such that the vertex $n-1$ adjacent to the vertex $n$, the resultant
graph is still denoted by $G$; applying Lemma \ref{a1} to $G$ such
that the vertex $n-2$ adjacent to the vertex $n-1$, and so on. The
graph $G(p,q;n)$ can be obtained. Applying Lemma \ref{a1} again,
each step above, the origin graph and its resultant graph have the
same determinant,  the result thus follows. \hfill
$\blacksquare$\vspace{3mm}

%
%
%
%
%
%
From Theorem \ref{a4}, to compute the determinant of the distance
matrix of a  bicyclic graph of order $n+p+q-1$ which contains
$\infty(p,k,q)$  as an induced subgraph, it is sufficient to compute
the determinant of the distance matrix of the graph $G(p,q;n)$. For
convenience, in the following denote by $D_n$ the distance matrix of
the graph $G(p,q;n)$.

\begin{theorem} \label{a5} Fixed the integers $p$ and $q$. If
 $n\ge 2$, then
 $$det D_n=-4det D_{n-1}-4det
D_{n-2}.$$
\end{theorem}
{\bf Proof.} Let the vertex set of $G(p,q;n)$ be
$\{1,2,\cdots,p+q+n-1\}$. Then $G(p,q;n-1)$ can be considered as the
induced subgraph of $G(p,q;n)$ by deleting the pendent vertex
$p+q+n-1$ and $G(p,q;n-2)$ can be considered as the induced subgraph
of $G(p,q;n)$ by deleting the pendent vertex $p+q+n-1$ together with
its neighbor $p+q+n-2$. Hence, $D_n$ can be partitioned as{\small
$$D_n=\left(\begin{array}{cccc}
D_{n-3}& d^T&d^T+\textbf{1}^T& d^T+2\textbf{1}^T\\
d& 0&1& 2\\
d+\textbf{1}&1&0& 1 \\
d+2\textbf{1}&2&1& 0
\end{array}\right),$$}where $(d~~0)$  denotes the row vector of
$D_{n-1}$ corresponding to the vertex $p+q+n-3$. Hence,
{\small \begin{eqnarray*}det D_n&=&\left|\begin{array}{cccc}
D_{n-3}& d^T&d^T+\textbf{1}^T& d^T+2\textbf{1}^T\\
d& 0&1& 2\\
d+\textbf{1}&1&0& 1 \\
d+2\textbf{1}&2&1& 0
\end{array}\right|\\&=&\left|\begin{array}{cccc}
D_{n-3}& d^T&\textbf{1}^T&\textbf{1}^T\\
d& 0&1& 1\\
d+\textbf{1}&1&-1& 1 \\
d+2\textbf{1}&2&-1& -1
\end{array}\right|=\left|\begin{array}{cccc}
D_{n-3}& d^T&\textbf{1}^T&\textbf{1}^T\\
d& 0&1& 1\\
\textbf{1}&1&-2& 0 \\
\textbf{1}&1&0& -2
\end{array}\right|\\&=&\left|\begin{array}{cccc}
 D_{n-3}& d^T&\textbf{1}^T&0\\
d& 0&1& 0\\
\textbf{1}&1&-2& 2 \\
0&0&2& -4
\end{array}\right|\\&=&-4\left|\begin{array}{ccc}
 D_{n-3}& d^T&\textbf{1}^T\\
d& 0&1\\
\textbf{1}&1&-2
\end{array}\right|-4\left|\begin{array}{cc}
 D_{n-3}& d^T\\
d& 0
\end{array}\right|\\&=&-4\left|\begin{array}{ccc}
 D_{n-3}& d^T&d^T+\textbf{1}^T\\
d& 0&1\\
d+\textbf{1}&1&0
\end{array}\right|-4\left|\begin{array}{cc}
 D_{n-3}& d^T\\
d& 0
\end{array}\right|\\&=&-4det D_{n-1}-4det D_{n-2}.\end{eqnarray*}}
The result follows. \hfill $\blacksquare$\vspace{3mm}

%
%

\begin{theorem} \label{a6} Fixed the integers $p$ and $q$.
Then $$det D_0=det D_1=0$$ if one of the integers $p$ and $q$ is
even; and \begin{eqnarray*} det D_0&=&\frac{(pq-1)(p+q)}{4}\\det
D_1&=&-\frac{1}{2}(p+q)(pq-1)-pq
\end{eqnarray*}
otherwise.
\end{theorem}

{\bf Proof.} Without loss of generality, suppose that, in
$G(p,q;0)$, $\{1,2,\cdots,p\}$ and $\{1,p+1,\cdots,p+q-1\}$ are
respectively the natural sequences of the vertex sets of the cycles
$C_p$ and $C_q$, and  in $G(p,q;1)$ the unique pendent vertex is
labeled as $p+q$. Then $D_1$ has the form as (2.1) and $D_0$ is the
submatrix of (2.1) by deleting the last row and the last column,
where $D^*$ and $D^{**}$ are respectively defined as Lemma \ref{a3}.
Hence from the proof of Lemma \ref{a3}, we have {\small
\begin{eqnarray*} det D_1&=&\left|\begin{array}{cccc}
0&d_1&d_{p+q-1}& 1\\
d_1^T& D^* &d_{p+q-1}\textbf{1}^T+\textbf{1}d_1^T& d_1^T+\textbf{1}^T\\
d_{p+q-1}^T& \textbf{1}d_{p+q-1}^T+d_1\textbf{1}^T &D^{**}& d_{p+q-1}^T+\textbf{1}^T\\
1& d_1+\textbf{1} & d_{p+q-1}+\textbf{1}&0
\end{array}\right|\\&=&\left|\begin{array}{cccc}
0&d_1&d_{p+q-1}& 1\\
d_1^T& D^*-\textbf{1}d_1^T-d_1\textbf{1}^T &0& d_1^T\\
d_{p+q-1}^T& 0 &D^{**}-\textbf{1}d_{p+q-1}^T-d_n\textbf{1}^T& d_{p+q-1}^T\\
1& d_1 & d_{p+q-1}&0
\end{array}\right|\\&=&\left|\begin{array}{cccc}
0&d_1&d_{p+q-1}& 1\\
d_1^T& D^*-\textbf{1}d_1^T-d_1\textbf{1}^T &0& 0\\
d_{p+q-1}^T& 0 &D^{**}-\textbf{1}d_{p+q-1}^T-d_n\textbf{1}^T& 0\\
1& 0& 0&-2
\end{array}\right|,
\end{eqnarray*}}
and {\small \begin{eqnarray*} det D_0&=&\left|\begin{array}{cccc}
0&d_1&d_{p+q-1}\\
d_1^T& D^*-\textbf{1}d_1^T-d_1\textbf{1}^T &0\\
d_{p+q-1}^T& 0 &D^{**}-\textbf{1}d_{p+q-1}^T-d_{p+q-1}\textbf{1}^T
\end{array}\right|.
\end{eqnarray*} }  We first
  discuss the matrix
$D^*-d_1\textbf{1}^T-\textbf{1}d_1^T$ and denote it by
$D^p=(d_{i,j})$ for simplify. Recall that we set $\{1,2,\cdots,p\}$
is the natural sequences of the vertex sets of $C_p$. Then for
$D^*=(d^*_{i,j})$ we have $d_{ij}^*=min\{p-|i-j|,|i-j|\}$, and
$$d_1=(1,2,\cdots,k-1,k,k-1,\cdots,2,1)$$ if $p=2k$;
$$d_1=(1,2,\cdots,k-1,k,k,k-1,\cdots,2,1)$$ if $p=2k+1$. Hence, for
$D^p=(d^p_{i,j})$, we have
$$ d^p_{i,j}=min\{p-|i-j|,|i-j|\}-min\{p-i+1,i-1\}-min\{p-j+1,j-1\}.\eqno{(3.1)}$$
Thus, if $p=2k$, {\small
$$D^p=\left(\begin{array}{cccccc|c|cccccc}
-2&-2&-2&\cdots&-2&-2&-2&0&0&0&\cdots&0&0\\
-2&-4&-4&\cdots&-4&-4&-4&-2&0&0&\cdots&0&0\\
\cdots&\cdots&\cdots&\cdots&\cdots&\cdots&\cdots&\cdots&\cdots&\cdots&\cdots&\cdots\\
-2&-4&-6&\cdots&4-2k&4-2k&4-2k&4-2k&6-2k&8-2k&\cdots&0&0\\
-2&-4&-6&\cdots&4-2k&2-2k&2-2k&2-2k&4-2k&6-2k&\cdots&-2&0\\
\hline
-2&-4&-6&\cdots&4-2k&2-2k&-2k&2-2k&4-2k&6-2k&\cdots&-4&-2\\
\hline
0&-2&-4&\cdots&4-2k&2-2k&2-2k&2-2k&4-2k&6-2k&\cdots&-4&-2\\
0&0&-2&\cdots&4-2k&4-2k&4-2k&4-2k&6-2k&8-2k&\cdots&-4&-2\\
\cdots&\cdots&\cdots&\cdots&\cdots&\cdots&\cdots&\cdots&\cdots&\cdots&\cdots&\cdots\\
0&0&0&\cdots&0&-2&-4&-4&-4&-4&\cdots&-4&-2\\
0&0&0&\cdots&0&0&-2&-2&-2&-2&\cdots&-2&-2
\end{array}\right),$$} and if $p=2k+1$,
{\small $$D^p=\left(\begin{array}{ccccc|ccccc}
-2&-2&-2&\cdots&-2&-1&0&0&\cdots&0\\
-2&-4&-4&\cdots&-4&-3&-1&0&\cdots&0
\\
-2&-4&-6&\cdots&-6&-5&-3&-1&\cdots&0
\\
\cdots&\cdots&\cdots&\cdots&\cdots&\cdots&\cdots&\cdots&\cdots&\cdots
\\
-2&-4&-6&\cdots&-2k&1-2k&3-2k&5-2k&\cdots&-1\\ \hline
-1&\cdots&-3&-5&1-2k&-2k&\cdots&-6&-4&-2\\
\cdots&\cdots&\cdots&\cdots&\cdots&\cdots&\cdots&\cdots&\cdots&\cdots\\
0&\cdots&-1&-3&-5&-6&\cdots&-6&-4&-2
\\
0&\cdots&0&-1&-3&-4&\cdots&-4&-4&-2
\\
0&\cdots&0&0&-1&-2&\cdots&-2&-2&-2
\end{array}\right).$$}

For $p=2k$, we have {\small \begin{eqnarray*}
&&\left|\begin{array}{cc}
0&d_1\\
d_1^T&D^p
\end{array}\right| \\&=&\left|\begin{array}{c|cccccc|c|cccccc}
 0&1&2&3&\cdots&k-2&k-1&k&k-1&k-2&k-3&\cdots&2&1\\
 \hline
1&-2&-2&-2&\cdots&-2&-2&-2&0&0&0&\cdots&0&0\\
2&-2&-4&-4&\cdots&-4&-4&-4&-2&0&0&\cdots&0&0\\
\cdots&\cdots&\cdots&\cdots&\cdots&\cdots&\cdots&\cdots&\cdots&\cdots&\cdots&\cdots&\cdots\\
k-2&-2&-4&-6&\cdots&4-2k&4-2k&4-2k&4-2k&6-2k&8-2k&\cdots&0&0\\
k-1&-2&-4&-6&\cdots&4-2k&2-2k&2-2k&2-2k&4-2k&6-2k&\cdots&-2&0\\
\hline
k&-2&-4&-6&\cdots&4-2k&2-2k&-2k&2-2k&4-2k&6-2k&\cdots&-4&-2\\
\hline
k-1&0&-2&-4&\cdots&4-2k&2-2k&2-2k&2-2k&4-2k&6-2k&\cdots&-4&-2\\
k-2&0&0&-2&\cdots&4-2k&4-2k&4-2k&4-2k&6-2k&8-2k&\cdots&-4&-2\\
\cdots&\cdots&\cdots&\cdots&\cdots&\cdots&\cdots&\cdots&\cdots&\cdots&\cdots&\cdots&\cdots\\
2&0&0&0&\cdots&0&-2&-4&-4&-4&-4&\cdots&-4&-2\\
1&0&0&0&\cdots&0&0&-2&-2&-2&-2&\cdots&-2&-2
\end{array}\right| \end{eqnarray*}}
{\small \begin{eqnarray*}&=&\left|\begin{array}{c|cccccc|c|cccccc}
0& 1&2&3&\cdots&k-2&k-1&k&k-1&k-2&k-3&\cdots&2&1\\
\hline
1&-2&-2&-2&\cdots&-2&-2&-2&0&0&0&\cdots&0&0\\
1&0&-2&-2&\cdots&-2&-2&-2&-2&0&0&\cdots&0&0\\
\cdots&\cdots&\cdots&\cdots&\cdots&\cdots&\cdots&\cdots&\cdots&\cdots&\cdots&\cdots&\cdots\\
1&0&0&0&\cdots&-2&-2&-2&-2&-2&-2&\cdots&0&0\\
1&0&0&0&\cdots&0&-2&-2&-2&-2&-2&\cdots&-2&0\\
\hline
1&0&0&0&\cdots&0&0&-2&-2&-2&-2&\cdots&-2&-2\\
\hline
1&0&-2&-2&\cdots&-2&-2&-2&-2&0&0&\cdots&0&0\\
1&0&0&-2&\cdots&-2&-2&-2&-2&-2&0&\cdots&0&0\\
\cdots&\cdots&\cdots&\cdots&\cdots&\cdots&\cdots&\cdots&\cdots&\cdots&\cdots&\cdots&\cdots\\
1&0&0&0&\cdots&0&-2&-2&-2&-2&-2&\cdots&-2&0\\
1&0&0&0&\cdots&0&0&-2&-2&-2&-2&\cdots&-2&-2
\end{array}\right|\\&=&\left|\begin{array}{cccccc|c|cccccc|c}
0& 1&1&1&\cdots&1&1&1&1&1&1&\cdots&1&1\\
\hline
1&-2&0&0&\cdots&0&0&0&0&0&0&\cdots&0&0\\
1&0&-2&0&\cdots&0&0&0&-2&0&0&\cdots&0&0\\
\cdots&\cdots&\cdots&\cdots&\cdots&\cdots&\cdots&\cdots&\cdots&\cdots&\cdots&\cdots&\cdots\\
1&0&0&0&\cdots&-2&0&0&0&0&-2&\cdots&0&0\\
1&0&0&0&\cdots&0&-2&0&0&0&0&\cdots&-2&0\\
\hline
1&0&0&0&\cdots&0&0&-2&0&0&0&\cdots&0&-2\\
\hline
1&0&0&-2&\cdots&0&0&0&0&-2&0&\cdots&0&0\\
\cdots&\cdots&\cdots&\cdots&\cdots&\cdots&\cdots&\cdots&\cdots&\cdots&\cdots&\cdots&\cdots\\
1&0&0&0&\cdots&-2&0&0&0&0&-2&\cdots&0&0\\
1&0&0&0&\cdots&0&-2&0&0&0&0&\cdots&-2&0\\
1&0&0&0&\cdots&0&0&-2&0&0&0&\cdots&0&-2
\end{array}\right|.\end{eqnarray*}}
Note that all operation above disconcern the first row and the first
column, thus $det \tilde{D}=0$, as well as $det D=0$, if one of the
integers $p$ and $q$ is even.

We now consider the case that both of $p$ and $q$ are odd. For
$p=2k+1$, we have {\small \begin{eqnarray*}
&&\left|\begin{array}{cc}
0&d_1\\
d_1^T&D^p
\end{array}\right|\\&=&\left|\begin{array}{c|ccccc|ccccc}
0&1&2&\cdots&k-1&k&k&k-1&\cdots&2&1\\
\hline
1&-2&-2&-2&\cdots&-2&-1&0&0&\cdots&0\\
2&-2&-4&-4&\cdots&-4&-3&-1&0&\cdots&0
\\
3&-2&-4&-6&\cdots&-6&-5&-3&-1&\cdots&0
\\
\cdots&\cdots&\cdots&\cdots&\cdots&\cdots&\cdots&\cdots&\cdots&\cdots&\cdots
\\
k&-2&-4&-6&\cdots&-2k&1-2k&3-2k&5-2k&\cdots&-1\\ \hline
k&-1&\cdots&-3&-5&1-2k&-2k&\cdots&-6&-4&-2\\
\cdots&\cdots&\cdots&\cdots&\cdots&\cdots&\cdots&\cdots&\cdots&\cdots&\cdots\\
3&0&\cdots&-1&-3&-5&-6&\cdots&-6&-4&-2
\\
2&0&\cdots&0&-1&-3&-4&\cdots&-4&-4&-2
\\
1&0&\cdots&0&0&-1&-2&\cdots&-2&-2&-2
\end{array}\right|\\&=&\left|\begin{array}{c|ccccc|ccccc}
0&1&2&\cdots&k-1&k&k&k-1&\cdots&2&1\\
\hline
1&-2&-2&-2&\cdots&-2&-1&0&0&\cdots&0\\
1&0&-2&-2&\cdots&-2&-2&-1&0&\cdots&0
\\
1&0&0&-2&\cdots&-2&-2&-2&-1&\cdots&0
\\
\cdots&\cdots&\cdots&\cdots&\cdots&\cdots&\cdots&\cdots&\cdots&\cdots&\cdots
\\
1&0&0&0&\cdots&-2&-2&-2&-2&\cdots&-1\\ \hline
1&-1&\cdots&-2&-2&-2&-2&\cdots&0&0&0\\
\cdots&\cdots&\cdots&\cdots&\cdots&\cdots&\cdots&\cdots&\cdots&\cdots&\cdots\\
1&0&\cdots&-1&-2&-2&-2&\cdots&-2&0&0
\\
1&0&\cdots&0&-1&-2&-2&\cdots&-2&-2&0
\\
1&0&\cdots&0&0&-1&-2&\cdots&-2&-2&-2
\end{array}\right|\\&=&\left|\begin{array}{c|cccccc|cccccc}
0&1&1&\cdots&1&1&1&1&1&\cdots&1&1&1\\
\hline
1&-2&0&0&\cdots&0&0&-1&0&0&\cdots&0&0\\
1&0&-2&0&\cdots&0&0&-1&-1&0&\cdots&0&0
\\
1&0&0&-2&\cdots&0&0&0&-1&-1&\cdots&0&0
\\
\cdots&\cdots&\cdots&\cdots&\cdots&\cdots&\cdots&\cdots&\cdots&\cdots&\cdots&\cdots&\cdots
\\
1&0&0&0&\cdots&-2&0&0&0&0&\cdots&-1&0\\
1&0&0&0&\cdots&0&-2&0&0&0&\cdots&-1&-1\\ \hline
1&-1&-1&\cdots&0&0&0&-2&0&\cdots&0&0&0\\
1&0&-1&\cdots&0&0&0&0&-2&\cdots&0&0&0\\
\cdots&\cdots&\cdots&\cdots&\cdots&\cdots&\cdots&\cdots&\cdots&\cdots&\cdots&\cdots&\cdots\\
1&0&0&\cdots&-1&-1&0&0&0&\cdots&-2&0&0
\\
1&0&0&\cdots&0&-1&-1&0&0&\cdots&0&-2&0
\\
1&0&0&\cdots&0&0&-1&0&0&\cdots&0&0&-2
\end{array}\right|\end{eqnarray*}\begin{eqnarray*}&=&\left|\begin{array}{ccc}
0&\textbf{1}&\textbf{1}\\
\textbf{1}^T&-2I_k&B_k\\
\textbf{1}^T&B_k^T&-2I_k
\end{array}\right|,\end{eqnarray*}}
where $B_k$ is defined as Lemma \ref{a0}. Similarly, let $q=2h+1$.
Then
$$\left(\begin{array}{cc}
0&d_{p+q-1}\\
d_{p+q-1}^T&D^{**}-d_{p+q-1}\textbf{1}^T-\textbf{1}d_{p+q-1}^T
\end{array}\right)=\left(\begin{array}{ccc}
0&\textbf{1}&\textbf{1}\\
\textbf{1}^T&-2I_h&B_h\\
\textbf{1}^T&B_h^T&-2I_h
\end{array}\right).$$
Hence, $$ det D_0=\left|\begin{array}{cccccc}
0&\textbf{1}& \textbf{1}&\textbf{1}&\textbf{1}\\
\textbf{1}^T& -2I_k&B_k&0&0\\
\textbf{1}^T& B^T_k&-2I_k&0&0\\
\textbf{1}^T&0&0& -2I_h&B_h\\
\textbf{1}^T&0&0& B^T_h&-2I_h
\end{array}\right|$$ and
$$ det D_1=\left|\begin{array}{cccccc}
0&\textbf{1}& \textbf{1}&\textbf{1}&\textbf{1}&1\\
\textbf{1}^T& -2I_k&B_k&0&0&0\\
\textbf{1}^T& B^T_k&-2I_k&0&0&0\\
\textbf{1}^T&0&0& -2I_h&B_h&0\\
\textbf{1}^T&0&0& B^T_h&-2I_h&0\\
1& 0&0&0&0&-2\\
\end{array}\right|$$
Let $C_k=\frac{1}{2}B_kB_k^T-2I_{k}$,
$F_k=\frac{1}{2}\textbf{1}B_k^T+\textbf{1}$. Because
\begin{eqnarray*} &&\left(\begin{array}{cccccc}
1&0&\frac{1}{2}\textbf{1}&0&\frac{1}{2}\textbf{1}&\frac{1}{2}\\
0& I_k&\frac{1}{2}B_k&0&0&0\\
0& 0&I_k&0&0&0\\
0&0&0& I_h&\frac{1}{2}B_h&0\\
0&0&0& 0&I_h&0\\
0& 0&0&0&0&1
\end{array}\right)\left(\begin{array}{cccccc}
0&\textbf{1}& \textbf{1}&\textbf{1}&\textbf{1}&1\\
\textbf{1}^T& -2I_k&B_k&0&0&0\\
\textbf{1}^T& B^T_k&-2I_k&0&0&0\\
\textbf{1}^T&0&0& -2I_h&B_h&0\\
\textbf{1}^T&0&0& B^T_h&-2I_h&0\\
1& 0&0&0&0&-2
\end{array}\right)\\&=&\left(\begin{array}{cccccc}
\frac{1+k+h}{2}&F_k&0&F_h&0&0\\
F_k^T& C_k&0&0&0&0\\
\textbf{1}^T& B^T_k&-2I_k&0&0&0\\
F_h^T&0&0& C_h&0&0\\
\textbf{1}^T&0&0& B^T_h&-2I_h&0\\
1& 0&0&0&0&-2
\end{array}\right),
\end{eqnarray*}
\begin{eqnarray*}
&&\left(\begin{array}{ccccc}
1&0&\frac{1}{2}\textbf{1}&0&\frac{1}{2}\textbf{1}\\
0& I_k&\frac{1}{2}B_k&0&0\\
0& 0&I_k&0&0\\
0&0&0& I_h&\frac{1}{2}B_h\\
0&0&0& 0&I_h
\end{array}\right)\left(\begin{array}{ccccc}
0&\textbf{1}& \textbf{1}&\textbf{1}&\textbf{1}\\
\textbf{1}^T& -2I_k&B_k&0&0\\
\textbf{1}^T& B^T_k&-2I_k&0&0\\
\textbf{1}^T&0&0& -2I_h&B_h\\
\textbf{1}^T&0&0& B^T_h&-2I_h
\end{array}\right)\\&=&\left(\begin{array}{ccccc}
\frac{k+h}{2}&F_k&0&F_h&0\\
F_k^T& C_k&0&0&0\\
\textbf{1}^T& B^T_k&-2I_k&0&0\\
F_h^T&0&0& C_h&0\\
\textbf{1}^T&0&0& B^T_h&-2I_h
\end{array}\right)
\end{eqnarray*}
and
\begin{eqnarray*}
&&\left(\begin{array}{cccccc}
1&-F_kC_k^{-1}&0&-F_hC_h^{-1}&0&0\\
0& I_k&0&0&0&0\\
0& 0&I_k&0&0&0\\
0&0&0& I_h&0&0\\
0&0&0& 0&I_h&0\\
0& 0&0&0&0&1
\end{array}\right)\left(\begin{array}{cccccc}
\frac{1+k+h}{2}&F_k&0&F_h&0&0\\
F_k^T& C_k&0&0&0&0\\
\textbf{1}^T& B^T_k&-2I_k&0&0&0\\
F_h^T&0&0& C_h&0&0\\
\textbf{1}^T&0&0& B^T_h&-2I_h&0\\
1& 0&0&0&0&-2
\end{array}\right)\\&=&\left(\begin{array}{cccccc}
\frac{1+k+h}{2}-F_kC_k^{-1}F_k^T-F_hC_h^{-1}F_h^T&0&0&0&0&0\\
F_k^T& C_k&0&0&0&0\\
\textbf{1}^T& B^T_k&-2I_k&0&0&0\\
F_h^T&0&0& C_h&0&0\\
\textbf{1}^T&0&0& B^T_h&-2I_h&0\\
1& 0&0&0&0&-2
\end{array}\right),
\end{eqnarray*}
we have
$$det
D_0=(-2)^{k+h}(\frac{1}{2}k+\frac{1}{2}h-F_kC_k^{-1}F_k^T-F_hC_h^{-1}F_h^T)|C_k||C_h|$$
and $$det
D_1=-2(-2)^{k+h}(\frac{k+h+1}{2}-F_kC_k^{-1}F_k^T-F_hC_h^{-1}F_h^T)|C_k||C_h|.$$
From Lemma \ref{a0}, we have  $det C_k=\frac{(2k+1)(-1)^k }{2^k}$
and
 $F_kC_k^{-1}F_k^T=-\frac{k}{2(2k+1)}.$
Recall that $p=2k+1$ and $q=2h+1$, hence\begin{eqnarray*}det
D_0&=&\frac{(2k+1)(2h+1)}{2}[k+\frac{k}{2k+1}+h+\frac{h}{2h+1}]\\&=&k(k+1)(2h+1)+h(h+1)(2k+1)\\&=&
\frac{(pq-1)(p+q)}{4}\end{eqnarray*} and
\begin{eqnarray*}det
D_1&=&-2\frac{(2k+1)(2h+1)}{2}[k+\frac{k}{2k+1}+h+1+\frac{h}{2h+1}]
\\&=&-(2k+1)(2h+1)(k+h+1+\frac{k}{2k+1}+\frac{h}{2h+1})
\\&=&-(2k+1)(2h+1)(k+h+1)-k(2h+1)-h(2k+1)\\&=&-\frac{1}{2}(p+q)(pq-1)-pq.\end{eqnarray*}
The result thus follows.\vspace{3mm}

Putting  Theorem \ref{a6} into Lemma \ref{a00}, we have the main
result of this paper as follows.
\begin{theorem} \label{a7} Let $G$  be  an arbitrary  bicyclic graph of order $p+q-1+n$
 which  contains  $\infty(p,k,q)$
 as an induced subgraph with $n\ge k-1$. Denote by $D$ the distance matrix of $G$.
 Then $det D=0$ if one of   integers $p$ and
 $q$ is even, and $$det D=[\frac{(pq-1)(p+q)}{4}+\frac{n}{2}pq](-2)^n\eqno{(3.2)}$$ otherwise.
 \end{theorem}

 \noindent{\bf Remark.} Theorem \ref{a7} can be considered as a generalization
for Graham and Pollacks'  result on  the determinant of trees
\cite{gp} and Bapat, Kirkland and Neumanns' result on  the
determinant of  a unicyclic graph \cite{bk}. We   view one vertex as
a cycle with length $1$, then each vertex can be viewed as an
$\infty$-graph $\infty(1,1,1)$ and thus each tree contains
$\infty(1,1,1)$ as its induced subgraph. Consequently, the distance
matrix of each tree of order $n$ has the same determinant of that of
the graph $G(1,1;n-1)$,   then from (3.2)
$$det D=det D(G(1,1;n-1))=\frac{n-1}{2}(-2)^{n-1}=-(n-1)(-2)^{n-2}$$
 which is coincide with the formula for the determinant of a tree. Let $G$ be a unicyclic graph of order $n+p$ whose
unique cycle has length $p$, then  such a unique cycle can be viewed
as an $\infty$-graph $\infty(p,1,1)$. Thus  the distance matrix of
such a graph has the same determinant of that of the graph
$G(p,1;n)$ by Theorem \ref{a7}. Hence, $det D=det D(G)=0$ if $p$ is
even, and $$det D=det
D(G)=[\frac{(p-1)(p+1)}{4}+\frac{n}{2}p](-2)^n$$
 if $p$ is odd from (3.2),
 which is coincide with
Theorems 3.4 and 3.7 of \cite{bk}.

\end{document}